

\input amstex

\documentstyle{amsppt}

\loadbold

\magnification=\magstep1

\pageheight{9.0truein}
\pagewidth{6.5truein}


\def\n{\noindent}

\def\s{\smallskip}
\def\h{\hskip .2truein}
\def\nxn{n{\times}n}

\def\nxn{n{\times}n}
\def\sgn{\operatorname{sgn}}
\def\trace{\operatorname{trace}}
\def\u{\bold{u}}
\def\v{\bold{v}}
\def\x{\bold{x}}
\def\w{\bold{w}}
\def\y{\bold{y}}
\def\z{\bold{z}}
\def\a{\bold{a}}
\def\b{\bold{b}}
\def\ZZ{\Bbb{Z}}
\def\QQ{\Bbb{Q}}
\def\Mlm{M_{\ell{\times}m}}
\def\Mlpm{M_{\ell + m}}
\def\Elm{E_{(\ell,m)}}
\def\Tlm{T_{(\ell,m)}}
\def\P{\Cal{P}}
\def\supd{\text{{\it SupDiag}}}
\def\subd{\text{{\it SubDiag}}}
\def\U{\bold{U}}
\def\T{\bold{T}}
\def\V{\bold{V}}
\def\W{\bold{W}}

\def\Y{\bold{Y}}

\def\D{\Cal{D}}
\def\Art{1}
\def\ArtSch{2}
\def\Artetal{3}
\def\BaySti{4}
\def\Bluetal{5}
\def\BroDuC{6}
\def\BroGoo{7}
\def\Cohetal{8}
\def\Dic{9}
\def\Jat{10}
\def\Kol{11}
\def\Letone{12}
\def\Lettwo{13}
\def\Pap{14}
\def\Pro{15}
\def\Ron{16}
\def\Row{17}

\topmatter

\title Effective Detection of Nonsplit Module Extensions \endtitle


\rightheadtext{Detection of Nonsplit Extensions}

\author Edward S. Letzter \endauthor

\abstract Let $n$ be a positive integer, and let $R$ be a finitely
presented (but not necessarily finite dimensional) associative algebra
over a computable field. We examine algorithmic tests for deciding (1)
if every at-most-$n$-dimensional representation of $R$ is semisimple,
and (2) if there exist nonsplit extensions of non-isomorphic
irreducible $R$-modules whose dimensions sum to no greater than
$n$. \endabstract

\address Department of Mathematics, Temple University, Philadelphia, PA 19122
\endaddress

\email letzter\@math.temple.edu \endemail

\thanks The author's research was supported in part by NSF grants
DMS-9970413 and DMS-0196236.  \endthanks

\endtopmatter

\document

\baselineskip = 12pt plus 2pt

\lineskip = 2pt minus 1pt

\lineskiplimit = 2pt

\head 1. Introduction \endhead 

If $R = k\{X_1,\ldots,X_s\}/\langle f_1,\ldots,f_t \rangle$ is a
finitely presented algebra over a field $k$, then it is easy to see
that the $n$-dimensional representations of $R$ amount to solutions to
a system of $tn^2$ commutative polynomial equations in $sn^2$
variables. Moreover, the $n$-dimensional irreducible representations
of $R$ can also be explicitly parametrized by finite systems of
commutative polynomial equations (cf\. \cite{\Art,
\Pro}). Consequently, the techniques of compuational algebraic
geometry (and in particular, Groebner basis methods) can be used to
study the $n$-dimensional representation theory of $R$
(cf\. \cite{\Letone, \Lettwo}); for example, the question of whether
or not $R$ has an irreducible $n$-dimensional representation can be
algorithmically decided (when $k$ is computable). In this paper we
consider algorithmic approaches to another fundamental question in the
representation theory of $R$: Do there exist nonsplit extensions of
finite dimensional $R$-modules?

We present effective procedures for deciding (1) if every
at-most-$n$-dimensional representation of $R$ is semisimple (i.e., if
there exist no nonsplit extensions of modules whose dimensions sum to
no greater than $n$), and (2) if there exists a nonsplit extension of
an $m$-dimensional irreducible representation of $R$ by a
non-isomorphic $\ell$-dimensional irreducible representation, for some
$\ell + m \leq n$. These procedures are indirect -- they do not give
the exact dimensions of the detected nonsplit extensions. However,
precise (and more costly) algorithms can be subsequently derived.

Our basic strategy is to reduce each of the considered representation
theoretic decision problems to the problem of deciding whether a particular
finite set of commutative polynomials has a common zero. Standard methods of
computational algebraic geometry can then be applied (in principle). A brief
discussion of the complexity of this approach is given in (2.6). The case
when $n = 2$, discussed in \S 5, provides an elementary illustration.

When $R$ is known beforehand to be finite dimensional over $k$, effective
methods for determining a linear basis for the Jacobson radical of $R$ have
been given in \cite{\Cohetal; \Dic; \Ron}.

\subhead Acknowledgement \endsubhead My thanks to the referee for
suggestions on clarifying the exposition.

\head 2. Preliminaries \endhead

In this section we develop our notation (which will remain fixed for
the remainder) and quickly review some necessary background.

\subhead 2.1 \endsubhead We assume throughout this paper that $\ell$,
$m$, and $n$ are positive integers, that $k$ is a field, that $K$ is a
field extension of $k$, that $f_1,\ldots, f_t$ are noncommutative
polynomials in the free associative $k$-algebra $k\{ X_1,\ldots,
X_s\}$, and that $R$ is the quotient algebra
$$k\{X_1,\ldots, X_s\}/\langle f_1,\ldots , f_t \rangle .$$
Let $d$ denote the maximum of the total degrees of the
$f_1,\ldots,f_t$.

\subhead 2.2 \endsubhead (i) We will use the term {\sl
indeterminate\/} only in reference to a variable in an (often tacitly
given) commutative polynomial ring. Unless otherwise designated, {\sl
polynomial\/} will refer only to a commutative polynomial.

(ii) Let $A$ be a $k$-algebra (algebras, modules, and homomorphisms
will always be assumed to be unital). If $a_1,\ldots,a_q \in A$, we
use $k\{a_1,\ldots,a_q\}$ to denote the $k$-subalgebra generated by
$a_1,\ldots,a_q$. 

Recall that every $K$-algebra automorphism $\tau$ of $M_n(K)$ is
{\sl inner\/} (i.e., there exists an invertible matrix $Q$ in $M_n(K)$
such that $\tau(a) = QaQ^{-1}$ for all $a \in M_n(K)$).

(iii) We let $M_n(K)$ denote the ring of $\nxn$ matrices with entries in
$K$, and we let $\Mlm(K)$ denote the $M_\ell(K)$-$M_m(K)$-bimodule of
$\ell{\times}m$ matrices. We identify $K^n$ with the left
$M_n(K)$-module of $n{\times}1$ matrices with entries in $K$.

Let $I_n$ denote the $\nxn$ identity matrix.  When $\ell < n$, we identify
$I_\ell$ with the $\nxn$ matrix $\bmatrix I_\ell & 0 \\ 0 & 0
\endbmatrix$. Let $\supd_n$ denote the $\nxn$ matrix with $1$'s on the
super-diagonal and $0$'s elsewhere, and let $\subd_n$ denote the transpose of
$\supd_n$. It is easy to verify that $\supd_n$ and $\subd_n$ generate $M_n(K)$
as a $K$-algebra.

(iv) We will use the expression {\sl ($n$-dimensional) representation
of $A$\/} only to refer to $k$-algebra homomorphisms $\rho\colon A
\rightarrow M_n(K)$; the representation is {\sl irreducible\/} when
$K\rho(A) = M_n(K)$.

This approach allows us to consider the $K$-representation theory of
$A$ while restricting our calculations to $k$; in our algorithmic
procedures below we will assume that $k$ is computable and that $K$ is
the algebraic closure of $k$. (Recall, if $K$ is the algebraic closure
of $k$, that a representation $\rho\colon R \rightarrow M_n(K)$ is
irreducible -- in the preceding sense -- if and only if the only
$K\rho(R)$-invariant subspaces of $K^n$ are $0$ and $K^n$ itself.)

(v) Two representations $\rho,\rho'\colon A \rightarrow M_n(K)$ are
{\sl equivalent\/} (or {\sl isomorphic\/}) provided there exists an
invertible matrix $Q \in M_n(K)$ such that $\rho'(a) = Q^{-1}\rho(a)Q$
for all $a \in A$. We will say that a representation $\rho$ of $A$ is
{\sl semisimple\/} if $K\rho(A)$ is semisimple as a $K$-algebra.

\subhead 2.3 \endsubhead (i) For $1 \leq \mu \leq s$, let $\x_\mu$ denote
the generic $\nxn$ matrix $(x_{ij}(\mu))$ (i.e., the $\nxn$ matrix
whose $ij$th entry is the indeterminate $x_{ij}(\mu)$), and set $\x =
(\x_1,\ldots,\x_s)$. Note that $R$ has an $n$-dimensional
representation if and only if the entries of $f_1(\x), \ldots ,
f_t(\x)$ have a common zero.

(ii) (Assume that $k$ is computable and that $K$ is the algebraic closure of
$k$.) Using standard techniques of computational commutative algebra, we can
check if $f_1(\x),\ldots,f_t(\x)$ have a common zero, and thereby decide
whether or not $R$ has an $n$-dimensional representation. Also, we can always
slightly simplify the computations by replacing one of the generic matrices
$(x_{ij}(\mu))$ with an upper triangular matrix (i.e., by setting
$x_{ij}(\mu)=0$ for $i > j$). Therefore, this procedure involves $tn^2$
polynomials, of degree at most $d$, in $sn^2 - (n^2 -n)/2$ variables. (Of
course, the specific relations defining $R$ may allow for further reductions.)

In all of the tests discussed below, we will assume that one of the generic
matrices has been similarly replaced with a generic upper triangular matrix.

\subhead 2.4 \endsubhead (i) Let $\P(n)$ denote the minimum
positive integer with the following property: For all positive
integers $q$, and for all $a_1,\ldots,a_q \in M_n(K)$, the $K$-algebra
$K\{a_1,\ldots,a_q\}$ is $K$-linearly spanned by products of the
$a_1,\ldots,a_q$ having length no greater than $\P(n)$. (The identity
matrix is a product of length zero.)

It is easy to check that $\P(n) \leq n^2 -1$, and in \cite{\Pap} it is
proved that $\P(n)$ is bounded above by a function in
$\Cal{O}(n^{3/2})$.

(ii) Let $\rho\colon R \rightarrow M_n(K)$ be a representation, and
set $\Lambda = K\rho(R)$. It follows from (i) that $\Lambda$ is
$K$-linearly spanned by the images of the monomials (in the $X_i$)
having length no greater than $\P(n)$. Also, the Cayley-Hamilton Theorem
tells us that the $n$th power of a matrix in $M_n(K)$ is a $K$-linear
combination of its lower powers. Therefore, $\Lambda$ is $K$-linearly
spanned by the image under $\rho$ of
$$\big\{ Y^{i_1}_1 \cdots Y^{i_p}_p  :  Y_1,\ldots,Y_p \in
\{1,X_1,\ldots,X_s\};  i_1 + \cdots + i_p \leq \P(n); 
i_1,\ldots,i_p < n  \big\} .$$

\subhead 2.5 \endsubhead For later comparison, we briefly mention two
algorithmic tests for detecting irreducible $n$-dimensional
representations. Let $\W=$
$$\big\{ \w^{i_1}_1\cdots \w^{i_p}_p : \w_1,\ldots,\w_p \in
\{I_n,\x_1,\ldots,\x_s\}; i_1 + \cdots + i_p \leq \P(n);
i_1,\ldots,i_p < n \big\} .$$
Assume (for the rest of this subsection) that $k$ is computable and
that $K$ is the algebraic closure of $k$.

(i) (Naive Irreducibility Test) For each choice of $\w_1,\ldots, \w_{n^2} \in
\W$ we can construct a subtest that returns ``true' if the entries of
$$f_1(\x),\ldots,f_t(\x), \quad y_1\w_1 + \cdots + y_{n^2}\w_{n^2} -
\supd_n, \quad z_1\w_1 + \cdots + z_{n^2}\w_{n^2} - \subd_n,$$
have a common zero, for indeterminates $y_i$ and $z_i$. The subtest
returns ``false'' if no common zero exists.

It follows immediately that the following are equivalent: (1) at least one of
the possible choices of $\w_1,\ldots,\w_{n^2}$ produces a ``true'' in the
subtest, (2) there exists an irreducible representation $R \rightarrow
M_n(K)$. (Of course, $\supd_n$ and $\subd_n$ can be replaced with any pair of
matrices in $M_n(k)$ that generate $M_n(K)$ as a $K$-algebra.)

Note that each subtest involves $(t+2)n^2$ polynomials in $(s+2)n^2 - (n^2
-n)/2$ variables. The degrees of $2n^2$ of these polynomials will be bounded by
$\P(n) + 1$, and the remaining degrees will be bounded by $d$.

(ii) Recall the {\sl $\nu$th standard identity\/},
$$s_\nu \quad = \quad \sum _{\sigma \in S_\nu} (\sgn \sigma)
Y_{\sigma(1)} \cdots Y_{\sigma(\nu)} \quad \in \quad
\ZZ\{Y_1,\ldots,Y_\nu\} .$$
(See, e.g., \cite{\Row}.) Observe that $s_\nu$ is multilinear and
alternating. Choose $\w_0,\ldots,\w_{2(n-1)} \in W$, and let $w$ be an
indeterminate. Consider a test that returns ``true'' if
$$w\trace(\w_0s_{2(n-1)}(\w_1,\ldots,\w_{2(n-1)})) - 1 $$
and the entries of $f_1(\x),\ldots,f_t(\x)$ have a common zero (and
returns ``false'' otherwise). In \cite{\Letone} it is shown that $R$
has an irreducible $n$-dimensional representation if and only if at
least one of these tests returns a ``true.''

Each subtest in this procedure will involve $tn^2 + 1$ polynomials in
$sn^2 - (n^2 -n)/2 + 1$ variables. One of these polynomials will have
degree $\P(n)^{2n-1} + 1$, and the remaining degrees will be bounded
by $d$.

\subhead 2.6 \endsubhead Koll\'ar's Sharp Effective Nullstellensatz
\cite{\Kol} offers a rough method, as follows, to compare the complexities of
the algorithms we encounter (cf\. \cite{\BaySti}).

(i) Let $q,r,d_1 \leq \cdots \leq d_q$ be positive integers, with no
$d_i = 2$, and with $r > 1$. Set
$$\D = \cases d_1\cdots d_q & q \leq r \\ d_1\cdots d_{r-1}d_q &
1 < r < q . \endcases $$

(ii) Let $g_1,\ldots,g_q \in k[x_1,\ldots,x_r]$, and suppose that $d_i =
\deg(g_i)$ for $1 \leq i \leq q$. In \cite{\Kol} it is shown that
$g_1,\ldots,g_q$ have no common zero (over the algebraic closure of $k$) if
and only if there exist $h_1,\ldots,h_q \in k[x_1,\ldots,x_r]$ such that
$h_1g_1 + \cdots + h_qg_q = 1$ and such that the degrees of the $g_ih_i$ are
no greater than $\D$. It is further shown in \cite{\Kol}, for arbitrarily
chosen $g_1,\ldots,g_q$ satisfying the given criteria, that this degree bound
is as small as possible.

(iii) Following \cite{\BaySti, \S 3} (cf\. \cite{\Bluetal, 1.2.5}), we
use $\D$ as a relative measure of the complexity of determining
whether $g_1,\ldots,g_q$ have a common zero. (In measuring $\D$ for
the systems below, we will simply -- and simplistically -- assume that
the degree of a quadratic polynomial is replaced by a $3$ in the
appropriate calculation.)

(iv) Let $u$ denote the minimum of $s$ and $t$. For the test deciding
whether $R$ has an $n$-dimensional representation (2.3ii), $\D \leq
d^{un^2}$.

(v) For convenience, in comparing costs of algorithms we will assume
that $\P(n) \geq d$. 

(vi) For the first irreducibility test (2.5i), we see that $\D \leq
d^{un^2}(\P(n)+1)^{2n^2}$. For the second (2.5ii), we see that $\D
\leq d^{un^2}(\P(n)^{2n-1} + 1)$. 

(vii) Unfortunately, the degree bounds in (iv) and (vi) involve
factors no smaller than $n$ raised to a polynomial in $n$. The degree
bounds we will encounter in later sections behave similarly. However,
the calculation of $\D$, following \cite{\Kol}, does not take into
account the specific representation-theoretic sources of the
polynomials occurring. We therefore ask: What are the minimum degree
complexities of $n$-dimensional representation-theoretic decision
problems?

\head 3. Semisimplicity Test \endhead

Let $A$ denote a $k$-algebra.

\subhead 3.1 \endsubhead Set
$$\Elm(K) = \bmatrix M_\ell(K) & \Mlm(K) \\ 0 & M_m(K) \endbmatrix,$$
a $K$-subalgebra of $\Mlpm(K)$. 

The next result will form the foundation for our semisimplicity test. The
proof will follow immediately from (3.7).

\proclaim{3.2 Proposition} Every at-most-$n$-dimensional
representation of $A$ is semisimple if and only if $\supd_{\ell + m}
\not\in K\rho(A)$ for all representations $\rho\colon A \rightarrow
\Elm(K) \subset \Mlpm(K)$ such that $\ell + m \leq n$. \endproclaim

\subhead 3.3 \endsubhead We will need some more notation.

(i) Associated to $\Elm(K)$ are canonical $K$-algebra homomorphisms
$\pi_\ell\colon \Elm(K) \rightarrow M_\ell(K)$ and $\pi_m\colon \Elm(K)
\rightarrow M_m(K)$.

(ii) Viewing $K^{\ell + m}$ as a left $\Elm(K)$-module, identify
$K^\ell$ with the submodule comprised of those column vectors having
only zero entries below the $\ell$th position. Further identify $K^m$
with the $\Elm(K)$-module factor $K^{\ell + m}/K^\ell$.

(iii) Set 
$$\Tlm(K) = \bmatrix 0 & \Mlm(K) \\ 0 & 0
\endbmatrix,$$
the Jacobson radical of $\Elm(K)$. 

\subhead 3.4 \endsubhead For the remainder of this section, assume
that $\rho\colon A \rightarrow \Elm(K)$ is a representation, that
$\Lambda = K\rho(A)$, and that $J$ is the Jacobson radical of
$\Lambda$. Also, let $\tau$ be an inner $K$-algebra automorphism of
$\Mlpm(K)$ such that $\tau(\Elm(K)) \subseteq \Elm(K)$. Of course,
$\tau\rho$ will be a representation of $A$ equivalent to $\rho$.

\subhead 3.5 \endsubhead (i) If the compositions $\pi_\ell\rho$ and
$\pi_m\rho$ are both irreducible, we will say that $\rho$ is an {\sl
$(\ell, m)$-extension of irreducible representations}; we will further
say that $\rho$ is a {\sl self extension\/} when $\pi_m\rho$ and
$\pi_\ell\rho$ are equivalent representations (and so $\ell = m$).

(ii) An $(\ell , m)$-extension of irreducible representations {\sl
splits\/} if it is semisimple. It easily follows from standard results
that every at-most-$n$-dimensional representation of $A$ is semisimple
if and only if all $(\ell,m)$-extensions of irreducible
representations of $A$ split, for all choices of $\ell + m \leq n$.

\proclaim{3.6 Lemma} Assume that $\rho$ is a nonsplit
$(\ell,m)$-extension of irreducible representations.

{\rm (i)} $J = \Tlm(K)$. 

{\rm (ii)} Suppose that $\rho$ is a self extension. Then 
we can choose $\tau$ such that
$$\tau(\Lambda) = \left\{ \bmatrix a & b \\ 0 & a \endbmatrix  : 
a \in M_\ell(K), b \in \Mlm(K)  \right\} .$$

{\rm (iii)} Suppose that $\rho$ is not a self extension. Then $\Lambda
= \Elm(K)$.
\endproclaim

\demo{Proof} By considering the composition series $0 \subsetneq
K^\ell \subsetneq K^{\ell + m}$, we see that $J \subseteq \Tlm(K)$.

Being a nonzero $M_\ell(K)$-$M_m(K)$-bimodule, $J$ is a nonzero left
module over
$$M_\ell(K) \otimes _K (M_m(K))^{\text{op}} \quad \cong \quad
M_{\ell m}(K) .$$
Consequently, $\dim _K J \geq \ell m$, and so $J =
\Tlm(K)$. Part (i) follows. Parts (ii) and (iii) follow easily from (i).
\qed\enddemo

\proclaim{3.7 Lemma} {\rm (i)} Suppose that $\rho$ is semisimple.
Then $\supd_{\ell + m} \not\in \tau(\Lambda)$.

{\rm (ii)} Suppose that $\rho$ is an $(\ell,m)$-extension of
irreducible representations. Then $\rho$ does not split if and only if
$\supd_{\ell + m} \in \tau(\Lambda)$ for some choice of
$\tau$. \endproclaim

\demo{Proof} (i) The semisimplicity of $\rho$ implies that $\Lambda$
embeds into $M_\ell(K) \oplus M_m(K)$. Therefore, the maximum index
of nilpotence of elements in $\Lambda$ is less than $\ell + m$.

(ii) The ``only if'' statement follows from (3.6), and the ``if''
statement follows from (i).
\qed\enddemo

\subhead 3.8 \endsubhead The following notation will be used in
the procedures presented in (3.9), (3.10), and (4.2).

(i) For positive integers $\ell,m,r$, we will let $\b_r(\ell,m)$ denote the
$(\ell + m){\times}(\ell + m)$ matrix whose
$$\text{$ij$th entry}  =  \cases \text{the indeterminate
$x_{ij}(r)$ if $i \leq \ell$ or $j \geq m$,} \\ \text{$0$ otherwise}
. \endcases $$

(ii) For positive integers $\ell,m,s$, we will let $\U(\ell,m,s)$
denote the set of all products $\a^{i_1}_1 \cdots \a^{i_p}_p$ such
that
$$\multline \a_1,\ldots,\a_p \in \{I_{\ell +
m},\b_1(\ell,m),\ldots,\b_s(\ell,m)\}, \\ i_1 + \cdots + i_p \leq
\P(\ell + m), \quad \text{and} \quad i_1,\ldots,i_p < \ell + m
.\endmultline$$
Furthmore, temporarily letting $\U = \U(\ell,m,s)$, we will let
$\pi_\ell(\U)$ denote $\{ \pi_\ell(\u) : \u \in \U \}$ and $\pi_m(\U)$
denote $\{ \pi_m(\u) : \u \in \U \}$.

\subhead 3.9 Semisimplicity Test \endsubhead (Assume that $k$ is
computable and that $K$ is the algebraic closure of $k$.) We now
describe a test for deciding whether every at-most-$n$-dimensional
representation of $R$ is semisimple. That the procedure works as
stated follows directly from (3.2). Retain the notation of (3.8), and
let $x_1,x_2,\ldots$ be indeterminates.

\s\n\null\hrulefill\null\s

\n{\bf Input}: $f_1,\ldots,f_t \in k\{X_1,\ldots,X_s\}$, positive
integer $n$ 

\s\n{\bf Output}: ``all semisimple'' if every
at-most-$n$-dimensional representation 

\n of $k\{X_1,\ldots,X_s\}/\langle f_1,\ldots,f_t\rangle$ is
semisimple; ``not all semisimple'' otherwise

\s\n {\bf Begin} 

\s\n\h {\bf For} $1 \leq \ell < m \leq n$ {\bf do}:

\s\n\h $q := \ell^2 + \ell m + m^2$

\s\n\h $\V$ $:=$ set of subsets of $\U(\ell,m,s)$ having
cardinality $q$

\s\n\h $W := 0$

\s\n\h\h {\bf While} $\V \ne \emptyset$ and $W = 0$ {\bf do:}
 
\s\n\h\h Choose $\V_i = \{\u_1,\ldots,\u_q\} \in \V$ 

\s\n\h\h {\bf If} the entries of

\s\centerline{$x_1\u_1+\cdots+\x_q\u_q-\supd_{\ell+m}$}

\s\centerline{$f_1\big(\b_1(\ell,m),\ldots,\b_s(\ell,m)\big),\; \ldots, \;
f_t\big(\b_1(\ell,m),\ldots,\b_s(\ell,m)\big)$,}

\s\n\h\h have a common zero over $K$ {\bf then} $W := 1$

\s\n\h\h {\bf Else} $\V := \V \setminus \{\V_i\}$

\s\n\h\h {\bf End}

\s\n\h {\bf End}

\s\n {\bf If} $W = 0$ {\bf then} {\bf return} ``all semisimple''

\s\n {\bf Else} {\bf return} ``not all semisimple''

\s\n {\bf End} 

\n\null\hrulefill\null\s

\n Note that the subtest within the while loop involves $(t+1)q$
polynomials in $(s+1)q - (\ell^2-\ell)/2 - (m^2-m)/2$ variables. The
degrees of $q$ of these polynomials will be bounded by $\P(\ell + m) +
1$, and the remaining degrees will be bounded by $d$. Following (2.6),
$\D \leq d^{uq}(\P(\ell + m) + 1)^q$.

\subhead 3.10 Nonsplit $(\ell,m)$-Extension Test \endsubhead (Assume
that $k$ is computable and that $K$ is the algebraic closure of $k$.)
We now combine (3.7) with (2.5) to devise a procedure for deciding,
for fixed $\ell$ and $m$, whether $R$ has a nonsplit $(\ell,
m)$-extension of irreducible representations. Retain the notation of
(3.8), and let $v$, $w$, and $y_1,y_2,\ldots$ be
indeterminates. (Note: While the following algorithm works as stated,
it would be reasonable in general to first check for existence of
$\ell$-dimensional and $m$-dimensional irreducible representations,
following (2.5).)

\goodbreak\s\n\null\hrulefill\null\s

\s\n{\bf Input}: $f_1,\ldots,f_t \in k\{X_1,\ldots,X_s\}$, positive
integers $\ell$ and $m$

\s\n{\bf Output}: ``yes'' if there exists a nonsplit
$(\ell,m)$-extension of irreducible representations of 
$k\{X_1,\ldots,X_s\}/\langle f_1,\ldots,f_t\rangle$; ``no'' otherwise

\s\n {\bf Begin}

\s\n $q := \ell^2 + \ell m + m^2$

\s\n $\U := \U(\ell,m,s)$

\s\n $\V$ $:=$ set of subsets of $\pi_\ell(\U)$ having cardinality
$2(\ell-1)$

\s\n $\W$ $:=$ set of subsets of $\pi_m(\U)$ having cardinality
$2(m-1)$

\s\n $\Y$ $:=$ set of subsets of $\U$ having cardinality $q$

\s\n $\T := \U{\times}\V{\times}\U{\times}\W{\times}\Y$ 

\s\n $Z := 0$

\s\n\h {\bf While} $Z = 0$ and $\T \ne \emptyset$ {\bf do}:

\s\n\h Choose $\T_i = \big(\v_0, \{\v_1,\ldots,\v_{2(\ell
-1)}\},\w_0,\{\w_1,\ldots,\w_{2(m-1)}\}, \{\y_1,\ldots,\y_q\} \big)
\in \T$

\s\n\h {\bf If} the entries of

\s\centerline{$v\trace(\v_0s_{2(\ell-1)}(\v_1,\ldots,\v_{2(\ell-1)}))
- 1$,}

\s\centerline{$w\trace(\w_0s_{2(m-1)}(\w_1,\ldots,\w_{2(m-1)})) -
1$,}

\s\centerline{$y_1\y_1 + \cdots + y_q\y_q - \supd_{\ell + m}$}

\s\centerline{$f_1\big(\b_1(\ell,m),\ldots,\b_s(\ell,m)\big),\;
\ldots, \; f_t\big(\b_1(\ell,m),\ldots,\b_s(\ell,m)\big)$,}

\s\n\h have a common zero {\bf then} $Z := 1$ 

\s\n\h {\bf Else} 

\s\n\h $\T := \T \setminus \{\T_i\}$

\s\n\h {\bf End}

\s\n {\bf If} $Z = 1$ then {\bf return} ``yes'' 

\s\n {\bf Else} {\bf return} ``no''

\s\n {\bf End}

\n\null\hrulefill\null\s

\n The subtest within the while loop involves $(t+1)q+2$ polynomials
in $(s+1)q + 2 -(\ell^2-\ell)/2 - (m^2-m)/2$ variables. The degrees of
$q$ of the polynomials are bounded by $\P(\ell + m) + 1$, the degree
of one of the polynomials is bounded by $\P(\ell)^{2\ell -1}+1$, and
the degree of one of the polynomials is bounded by $\P(m)^{2m-1} +
1$. The remaining degrees are bounded by $d$. Following (2.6), $\D
\leq (\P(\ell)^{2\ell -1}+1)(\P(m)^{2m -1}+1)(\P(\ell +
m)+1)^qd^{uq}$.

\head 4. Nonsplit extensions of distinct irreducible representations
\endhead

\proclaim{4.1 Proposition} Let $A$ be a $k$-algebra. The following statements
are equivalent: {\rm (i)} There exists a nonsplit $(\ell, m)$-extension of
inequivalent irreducible representations of $A$ for some $\ell + m \leq n$.
{\rm (ii)} For some $\ell + m \leq n$, there exists a representation $\rho
\colon A \rightarrow \Elm(K)$ for which $\supd_{\ell + m},I_\ell \in
K\rho(A)$.
\endproclaim

\demo{Proof} (i)$\Rightarrow$(ii): Follows from (3.6iii).

(ii)$\Rightarrow$(i): Set $\Lambda = K\rho(A)$. If $K^{\ell + m}$ is
decomposable as a left $\Lambda$-module, then $\Lambda$ embeds into
$M_\mu(K) \oplus M_\nu(K)$, for some $\mu, \nu < \ell + m$, implying
that $\Lambda$ cannot contain an element whose index of nilpotence is
$\ell + m$. Therefore, since $\supd_{\ell + m} \in \Lambda$, we see
that $K^{\ell + m}$ is an indecomposable $\Lambda$-module.

Now let $M$ be the $\Lambda$-submodule $\Lambda I_\ell K^{\ell + m}$ of
$K^{\ell + m}$, and set $N = K^{\ell + m}/M$. Since $\Lambda$ is a subalgebra
of $\Elm(K)$, we see that both $M$ and $N$ are nonzero. It follows from the
preceding paragraph that the exact sequence $0 \rightarrow M \rightarrow
K^{\ell + m} \rightarrow N \rightarrow 0$ is a nonsplit extension of
$\Lambda$-modules. Therefore, there exists a nonsplit extension of $L'$ by $L$
for some simple $\Lambda$-module subfactor $L$ of $M$ and simple
$\Lambda$-module subfactor $L'$ of $N$. Note, however, that $I_\ell$ acts as
the identity on $L$ and that $I_\ell L' = 0$. Therefore, $L$ and $L'$ cannot
be isomorphic as $\Lambda$-modules.

Consequently, for some $1 \leq \ell' \leq \ell$ and $1
\leq m' \leq m$, there exists a nonsplit $(\ell',m')$-extension of
inequivalent irreducible representations $\rho'\colon A \rightarrow
E_{(\ell',m')}(K)$. \qed\enddemo

\subhead 4.2 Nonsplit Non-Self Extension Test \endsubhead (Assume that
$k$ is computable and that $K$ is the algebraic closure of $k$.)
Retain the notation of (3.8), and let $x_1,x_2,\ldots$ and
$y_1,y_2,\ldots$ be indeterminates. We can now describe a test, as
follows, for determining the existence of nonsplit extensions of
inequivalent irreducible representations. That the procedure works as
stated follows directly from (4.1).

\goodbreak\s\n\null\hrulefill\null\s

\n{\bf Input}: $f_1,\ldots,f_t \in k\{X_1,\ldots,X_s\}$, positive
integer $n$ 

\s\n{\bf Output}: ``yes'' if there exists a nonsplit
$(\ell,m)$-extension of inequivalent irreducible representations for
some $\ell + m \leq n$; ``no'' otherwise

\s\n {\bf Begin} 

\s\n\h {\bf For} $1 \leq \ell < m \leq n$ {\bf do}:

\s\n\h $q := \ell^2 + \ell m + m^2$

\s\n\h $\V$ $:=$ set of subsets of $\U(\ell,m,s)$ having
cardinality $q$

\s\n\h $W := 0$

\s\n\h\h {\bf While} $\V \ne \emptyset$ and $W = 0$ {\bf do:}
 
\s\n\h\h Choose $\V_i = \{\v_1,\ldots,\v_q\} \in \V$ 

\s\n\h\h {\bf If} the entries of

\s\centerline{$x_1\v_1+\cdots+x_q\v_q-\supd_{\ell+m}$}

\s\centerline{$y_1\v_1 + \cdots + y_q\v_q - I_\ell$}

\s\centerline{$f_1\big(\b_1(\ell,m),\ldots,\b_s(\ell,m)\big),\; \ldots, \;
f_t\big(\b_1(\ell,m),\ldots,\b_s(\ell,m)\big)$,}

\s\n\h\h have a common zero over $K$ {\bf then} $W := 1$

\s\n\h\h {\bf Else} $\V := \V \setminus \{\V_i\}$

\s\n\h\h {\bf End}

\s\n\h {\bf End}

\s\n {\bf If} $W = 1$ {\bf then} {\bf return} ``yes''

\s\n {\bf Else} {\bf return} ``no''

\s\n {\bf End} 

\n\null\hrulefill\null\s

\n The subtest within the while loop involves $(t+2)q$ polynomials in
$(s+2)q - (\ell^2-\ell)/2 - (m^2-m)/2$ variables. The degrees of $2q$
of these polynomials will be bounded by $\P(\ell + m) + 1$, and the
remaining degrees will be bounded by $d$. Following (2.6), $\D \leq
d^{uq}(\P(\ell + m) + 1)^{2q}$.

\subhead 4.3 \endsubhead We leave to the reader the construction of a
test that decides the existence of a nonsplit $(\ell,m)$-extension of
inequivalent irreducible representations, for fixed $\ell$ and $m$.

\head 5. Example: Nonsplit Extensions of One-Dimensional
Representations \endhead

As an elementary (and easy) illustration of the methods of the preceding
sections, we consider the case when $\ell = m = 1$. Nonsplit extensions of
one-dimensional representations play an important role in the study of many
natural classes of finitely presented algebras -- for example, in the study of
solvable Lie algebras (cf., e.g., \cite{\BroDuC}) and quantum function
algebras (e.g., \cite{\BroGoo}).

Assume that $k$ is computable and that $K$ is the algebraic closure of
$k$. Recall that $R = k\{X_1,\ldots,X_s\}/\langle f_1,\ldots,f_t
\rangle$.

\subhead 5.1 \endsubhead (i) For $1 \leq r \leq s$, set
$$\b_r = \bmatrix x_{11}(r) & x_{12}(r) \\ 0 & x_{22}(r)
\endbmatrix .$$

(ii) Following (2.5), and noting that $\P(2) \leq 3$, we set
$$\V = \big\{  I_n  \big\} \cup \big\{  \b_1,\ldots,\b_s  \big\}
\cup \big\{  \b_\alpha\b_\beta  :  \alpha \ne \beta  \big\}
\cup \big\{  \b_\alpha\b_\beta\b_\gamma  :  \alpha \ne \beta \ne
\gamma  \big\}.$$

(iii) Let $a_1$, $a_2$, and $a_3$ be indeterminates. By (3.9), there
exists a nonsplit extension of one-dimensional representations of $R$
if and only if the polynomial entries of
$$f_1(\b_1,\ldots,\b_s),\ldots,f_t(\b_1,\ldots,\b_s), \quad a_1\v_1 +
a_2\v_2 + a_3\v_3 - \bmatrix 0 & 1 \\ 0 & 0 \endbmatrix$$
have a common zero for some choice of distinct $\v_1,\v_2,\v_3 \in
\V$.

(iv) Let $a_1$, $a_2$, $a_3$, $b_1$, $b_2$, and $b_3$ be
indeterminates. By (4.2), there exists a nonsplit extension of
inequivalent one-dimensional representations of $R$ if and only if the
polynomial entries of
$$\multline f_1(\b_1,\ldots,\b_s),\ldots,f_t(\b_1,\ldots,\b_s), \\
\quad a_1\v_1 + a_2\v_2 + a_3\v_3 - \bmatrix 0 & 1 \\ 0 & 0
\endbmatrix, \\ \quad b_1\v_1 + b_2\v_2 + b_3\v_3 - \bmatrix 1 & 0 \\
0 & 0 \endbmatrix \endmultline$$
have a common zero for some choice of distinct $\v_1,\v_2,\v_3 \in \V$.

(v) If $s = 3$ then $\vert \V \vert = 22$, and there are ${22 \choose
3} = 1540$ cases to check in (iii) and (iv). It is not unusual for the
first interesting cases of a given class of algebras to require three
generators or fewer -- well-known occurrences of this phenomenon include
the enveloping algebra of $sl_2$, the enveloping algebra of the
Heisenberg Lie algebra, and the three-dimensional regular algebras of
\cite{\ArtSch; \Artetal}. 

\subhead 5.3 \endsubhead We conclude by considering two concrete
examples. All of the computations mentioned below were performed using
Macaulay2 on a personal computer (4 GB RAM). Let
$$\x = \bmatrix x_{11} & x_{12} \\ 0 & x_{22} \endbmatrix, \quad \y =
\bmatrix y_{11} & y_{12} \\ 0 & y_{22} \endbmatrix, \quad \z =
\bmatrix z_{11} & z_{12} \\ 0 & z_{22} \endbmatrix.$$

(i) Set
$$R = \QQ\{X,Y,Z\}/\langle XY-YX-Z,XZ-ZX,YZ-ZY \rangle,$$
the universal enveloping algebra of the (nilpotent) Heisenberg Lie
algebra. It follows from well known abstract arguments that $R$ does
not have nonsplit extensions of inequivalent one-dimensional
representations but does have nonsplit self extensions of
one-dimensional representations; see, for example, \cite{\Jat}.

Evaluating all 1540 cases, we were easily able to check that the
entries of
$$\multline \x\y-\y\x-\z, \quad \x\z-\z\x, \quad \y\z-\z\y \\ \quad
a_1\v_1 + a_2\v_2 + a_3\v_3 - \bmatrix 0 & 1 \\ 0 & 0 \endbmatrix, \\
\quad b_1\v_1 + b_2\v_2 + b_3\v_3 - \bmatrix 1 & 0 \\ 0 & 0
\endbmatrix, \endmultline$$
have no common zeros for indeterminates $a_1, a_2, a_3, b_1, b_2, b_3$
and all choices of $\v_1, \v_2, \v_3 \in \V$. We thus recovered the
fact that $R$ has no nonsplit extensions of distinct one-dimensional
representations. Next, evaluating all 1540 cases of
$$\x\y-\y\x-\z, \quad \x\z-\z\x, \quad \y\z-\z\y \quad a_1\v_1 +
a_2\v_2 + a_3\v_3 - \bmatrix 0 & 1 \\ 0 & 0 \endbmatrix, $$
we found that there exists a common zero -- indicating the presence of
a nonsplit self extension -- in 980 instances.

(ii) Now set
$$R = \QQ\{X,Y,Z\} / \langle XY-YX-Y, XZ-ZX, YZ-ZY\rangle,$$
an enveloping algebra of a solvable-but-not-nilpotent Lie algebra. Here it
follows from abstract considerations that $R$ has nonsplit self
and non-self extensions of one-dimensional representations (again see,
e.g., \cite{\Jat}).

Testing all 1540 cases of
$$\x\y - \y\x - \y, \quad \x\z-\z\x, \quad \y\z-\z\y, \quad a_1\v_1 +
a_2\v_2 + a_3\v_3 - \bmatrix 0 & 1 \\ 0 & 0 \endbmatrix,$$
for $\v_1,\v_2,\v_3 \in \V$, we found that there exists a common zero
-- indicating the presence of a nonsplit extension -- in 1539
instances. Only in the case $\{\v_1,\v_2,\v_3\} = \{I_2, \y\x\y,
\y\z\y\}$ did there not exist a common zero. Testing
$$\multline \x\y - \y\x - \y, \quad \x\z-\z\x, \quad \y\z-\z\y, \\
\quad a_1\v_1 + a_2\v_2 + a_3\v_3 - \bmatrix 0 & 1 \\ 0 & 0
\endbmatrix, \\ \quad b_1\v_1 + b_2\v_2 + b_3\v_3 - \bmatrix 1 & 0 \\
0 & 0 \endbmatrix, \endmultline$$
we found that there exists a common zero -- indicating the presence of
a nonsplit non-self extension -- in 650 instances.

\subhead 5.3 \endsubhead Unfortunately, at this time, we are unaware
of general methods for significantly simplifying the computations
involved in the procedures described in this paper. Systematic studies
of more practical approaches to these or related tests are left for
future work.

\enddocument